\documentclass[letter, 10pt, conference]{ieeeconf}      

\IEEEoverridecommandlockouts        

\usepackage{etex}
\usepackage[vlined,ruled]{algorithm2e}
\usepackage[dvipsnames]{xcolor}
\usepackage{pgfpages}
\usepackage[cmex10]{amsmath}		
\usepackage{amssymb, mathrsfs}
\usepackage{cite}
\usepackage{url}

\usepackage{dsfont}
\usepackage{siunitx}
\usepackage{verbatim}									
\usepackage{arydshln}
\usepackage{mathtools}									
\usepackage{siunitx}									
\usepackage{mdframed}

\usepackage{graphicx}
\usepackage{epstopdf}

\usepackage[font=small]{caption}
\usepackage{subcaption}
\usepackage[activate={true,nocompatibility},final,tracking=true,kerning=true,spacing=true,factor=1100,stretch=10,shrink=10]{microtype}

\usepackage{array}
\usepackage{multirow}
\usepackage{enumerate}
\usepackage{booktabs}

\graphicspath{{./fig/}}

\usepackage{pgfplots}
\pgfplotsset{width=8cm,compat=newest}
\newcommand*{\damping}{0.006}%
\newcommand*{\freq}{25}%
\pgfmathsetmacro{\freqd}{sqrt(1-(\damping)^2)*\freq}%

\pgfplotsset{
    standard/.style={
    axis x line=middle,
    axis y line=middle,
    enlarge x limits=0.15,
	enlarge y limits=0.15,
	every axis plot post/.style={mark options={fill=black}},
	}
}

\pgfplotsset{%
    ,compat=1.12
    ,every axis x label/.style={at={(current axis.right of origin)},anchor=north west}
    ,every axis y label/.style={at={(current axis.above origin)},anchor=north east}
    }

\usepackage{tikz}
\usetikzlibrary{arrows}
\usetikzlibrary{decorations.pathmorphing}
\usetikzlibrary{decorations.markings}
\usetikzlibrary{snakes}
\usetikzlibrary{matrix}
\usetikzlibrary{calc}
\usetikzlibrary{patterns}
\usetikzlibrary{positioning}
\usetikzlibrary{shapes}
\usetikzlibrary{shapes.symbols}
\usetikzlibrary{shapes.geometric}
\usetikzlibrary{shadows.blur}
\usetikzlibrary{shapes.arrows}
\usetikzlibrary{shapes.multipart}
\usetikzlibrary{shapes.callouts}
\usetikzlibrary{shapes.misc}

\tikzstyle{every node}=[font=\small]
\tikzstyle{every path}=[line width=0.8pt,line cap=round,line join=round]

\usepackage[american,smartlabels]{circuitikz}
\ctikzset{bipoles/thickness=1}
\ctikzset{bipoles/length=0.8cm}
\ctikzset{bipoles/diode/height=.375}
\ctikzset{bipoles/diode/width=.3}
\ctikzset{tripoles/thyristor/height=.8}
\ctikzset{tripoles/thyristor/width=1}
\ctikzset{bipoles/vsourceam/height/.initial=.7}
\ctikzset{bipoles/vsourceam/width/.initial=.7}

\newcommand{\real}{\mathbb{R}}

\newcommand{\setdef}[2]{\{#1 \;|\; #2\}}

\newcommand{\intersection}{\ensuremath{\operatorname{\cap}}}

\DeclareMathOperator*{\minimize}{minimize} 									

\newcommand{\vect}[1]{\mathbbold{#1}}
\newcommand{\vones}[1][]{\vect{1}_{#1}}

\DeclareSymbolFont{bbold}{U}{bbold}{m}{n}
\DeclareSymbolFontAlphabet{\mathbbold}{bbold}

\newcommand{\map}[3]{#1: #2 \rightarrow #3}

\renewcommand{\theenumi}{(\roman{enumi})}

\usepackage{soul}

\newcommand{\define}{\coloneqq}

\DeclareMathOperator{\subto}{subject~to}

\newcommand\oprocendsymbol{\hbox{$\square$}}
\newcommand\oprocend{\relax\ifmmode\else\unskip\hfill\fi\oprocendsymbol}

\newtheorem{theorem}{Theorem}[section]
\newtheorem{lemma}[theorem]{Lemma}

\newtheorem{proposition}[theorem]{Proposition}
\newtheorem{problem}{Problem}[section]

\newtheorem{assumption}{Assumption}[section]

\newenvironment{pfof}[1]{\vspace{1ex}\noindent{\itshape Proof of
    #1:}\hspace{0.5em}} {\hfill\oprocend\vspace{1ex}}

\newif\ifforstudents

\usepackage{makecell}

\usepackage{cuted}
\usepackage{arydshln}
\newcommand{\T}{\mathsf{T}}

\title{\bf Low-Gain Stabilizers for Linear-Convex Optimal Steady-State Control %
\thanks{
}}

\author{John W. Simpson-Porco%
\thanks{J.~W.~Simpson-Porco is with the Department of Electrical and Computer Engineering, University of Toronto, 10 King's College Road,
Toronto, ON, M5S 3G4, Canada. Email: {\tt jwsimpson@ece.utoronto.ca}. This work was supported in part by the NSERC Discovery Grant  RGPIN-2017-04008.}
}

\begin{document}
\maketitle
\thispagestyle{empty}
\pagestyle{empty}


\begin{abstract}
We consider the problem of designing a feedback controller which robustly regulates an LTI system to an optimal operating point in the presence of unmeasured disturbances. A general design framework based on so-called optimality models was previously put forward for this class of problems, effectively reducing the problem to that of stabilization of an associated nonlinear plant. This paper presents several simple and fully constructive stabilizer designs to accompany the optimality model designs from \cite{LSPL-JWSP-EM:18l}. The designs are based on a low-gain integral control approach, which enforces time-scale separation between the exponentially stable plant and the controller. We provide explicit formulas for controllers and gains, along with LMI-based methods for the computation of robust/optimal gains. The results are illustrated via an academic example and an application to power system frequency control. 
\end{abstract}

\section{Introduction}
\label{Section: Introduction}


Many practical engineering systems are subject to both dynamic response specifications and minimum-cost steady-state operation specifications. Achieving the former typically relies on both good process design and the design of one or more accompanying feedback controllers. In contrast, the latter optimality {criterion} is most often addressed in a hierarchical fashion, wherein actuator set-points are scheduled using a model and then forwarded to the lower-level feedback controllers as references. The quality of set-points computed via this procedure may however be poor, due to inaccuracies in the steady-state system model and due to discrepancies between forecasted and real-time disturbances. Such inaccuracies will lead to sub-optimal system operation, and may lead to violation of desired operating constraints; these issues motivate the incorporation of additional real-time feedback into the computation of optimal set-points.

Perhaps the most well-known framework which blends the control and optimization layers is model-predictive control (MPC), wherein simultaneous control and optimization is achieved via repeated solution of dynamic optimization problems \cite{JBR-DQM:09}. An alternative approach \textemdash{} tailored towards applications where low-level dynamic controllers are already in place \textemdash{} is to directly incorporate measurement feedback into the set-point scheduling process, such that the closed-loop system converges to a cost-minimizing operating point. By incorporating feedback, sensitivity to steady-state model uncertainty can be reduced, and constraint violation can be eliminated. While not as universally applicable as MPC, this approach has the benefit of producing simple and explicit controller designs, which are intuitive and interpretable to domain experts in particular applications (e.g., \cite{MF-CL-SL:13,AH-SB-GH-FD:16,ZT-EE-JWSP-EF-MP-HH:20l}).


The general problem described above been studied recently under different names, including dynamic KKT controllers \cite{AJ-ML-PPJvdB:09b}, feedback-based optimization \cite{MC-EDA-AB:18,SM-AH-SB-GH-FD:18,MC-JWSP-AB:19c,VH-AH-LO-SB-FD:21,GB-DLM-MHDB-SB-JL-FD:21}, autonomous optimization \cite{AH-SB-GH-FD:21,AH-SB-GH-FD:21b}, or \textemdash{} the term we will use here \textemdash{} optimal steady-state (OSS) control \cite{XZ-AP-NL:18,LSPL-ZEN-EM-JWSP:18e,LSPL-JWSP-EM:18l}; see also \cite{EDA-SVD-GBG:15,GB-MV-JC-ED:21,GB-JC-JIP-ED:21,JWSP:21a}. While the problem setups in the above references vary significantly in terms of the assumptions on the plant model and on the steady-state optimization problem, the general goal is the design of simple feedback laws which compute the set-points required to guide the plant towards a cost-minimizing and constraint-satisfying operating point.

\smallskip

\emph{Contributions:} This paper continues the development of the OSS control framework put forward in \cite{LSPL-ZEN-EM-JWSP:18e,LSPL-JWSP-EM:18l}, which will be reviewed in Section \ref{Sec:Main}. In brief, the framework in \cite{LSPL-JWSP-EM:18l} was shown to provide substantial flexibility to the designer, enabling customization of the control architecture (centralized vs. distributed) through judicious selection of several gain matrices. However, stabilizability guarantees were provided only for the case of quadratic steady-state optimization problems. Here, focusing on the practically-relevant case of a stable LTI plant, we provide stabilizer designs to accompany the ``optimality models'' from \cite{LSPL-JWSP-EM:18l} (Section \ref{Sec:ClosedLoop}). This renders the approach of \cite{LSPL-JWSP-EM:18l} fully constructive without the restriction of quadratic costs, {and} while retaining the design flexibility inherent {\cite{LSPL-JWSP-EM:18l}}. As a secondary contribution, the development of optimality models here (Section \ref{Sec:Optimality}) is more streamlined than in \cite{LSPL-JWSP-EM:18l}. Taken together, our results provide a library of constructive solutions for solving the OSS control problem in the framework of \cite{LSPL-JWSP-EM:18l} {for exponentially stable plants}.\footnote{Further results can be found in the extended version of this paper \cite{JWSP:22d-extended}.}





\section{Problem Formulation}
\label{Sec:Main}

Consider the finite-dimensional continuous-time LTI plant
\begin{equation}\label{Eq:LTI2}
\begin{aligned}
\dot{x} &= Ax + Bu + B_ww\\
z &= Cx + Du + D_ww
\end{aligned}
\end{equation}
with state $x(t) \in \real^n$, measured output $z(t) \in \real^r$, control input $u(t) \in \real^m$, and \emph{constant} exogenous signal $w \in \mathcal{W} \subset \real^{n_w}$, modelling reference signals and (potentially, unmeasured) external disturbances. As we will be subsequently pursuing a low-gain design philosophy \cite{EJD:76,JWSP:20a}, we assume that \eqref{Eq:LTI2} is internally exponentially stable, i.e., $A$ is Hurwitz. This stability may be inherent to the plant model, or may have been achieved through a previous compensator design. In this case, the steady-state input-output mapping of \eqref{Eq:LTI2} is
\begin{equation}\label{Eq:SteadyStateMap}
\bar{z} = G_u\bar{u} + G_ww
\end{equation}
{in which $\bar{z},\bar{u}$ denote equilibriums value of $z$ and $u$, and}
\[
\begin{aligned}
G_u &= -CA^{-1}B+D \in \real^{r \times m}\\
G_w &= -CA^{-1}B_w + D_w \in \real^{r \times n_w}
\end{aligned}
\]
denote the DC gain matrices of \eqref{Eq:LTI2}. Our controller design specification, which we refer to as the \emph{optimal steady-state} (OSS) control problem, roughly follows that in \cite{LSPL-JWSP-EM:18l}. 

\smallskip

\begin{problem}[\bf OSS Control]\label{Prob:OSS}
For the plant \eqref{Eq:LTI2}, design an output-feedback controller such that for each $w \in \mathcal{W}$, the closed-loop system possesses a unique locally exponentially stable equilibrium point, and such that the equilibrium pair $(\bar{u},\bar{z})$ is an optimal point of the problem
\begin{subequations}\label{Eq:OSS}
\begin{align}
\label{Eq:OSS-1}
&\minimize_{\bar{z},\bar{u}}& \quad &f_0(\bar{u}) + g_0(\bar{z})\\
\label{Eq:OSS-2}
& \subto &\quad &\bar{z} = G_u\bar{u} + G_ww\\
\label{Eq:OSS-3}
&&\quad &0 = H_{z}\bar{z} + H_{u}\bar{u} + H_ww.
\end{align}
\end{subequations}
\end{problem}
\medskip

In the optimization problem \eqref{Eq:OSS}, $\map{f_0}{\mathcal{U}}{\real}$ and $\map{g_0}{\mathcal{Z}}{\real}$ are convex objective functions to be minimized (e.g., operational costs, measures of tracking error), where $\mathcal{U} \subseteq \real^m$ and $\mathcal{Z} \subseteq \real^r$ are closed and non-empty convex sets. The functions $f_0, g_0$ may incorporate barrier or penalty terms for enforcement of inequality constraints. The constraint \eqref{Eq:OSS-2} is the steady-state constraint imposed by the dynamic system \eqref{Eq:LTI2}. Finally, \eqref{Eq:OSS-3} represents $n_{\rm c}$ additional affine \emph{engineering constraints} imposed by the designer, with $H_z \in \real^{n_{\rm c} \times r}$, $H_u \in \real^{n_{\rm c} \times m}$, and $H_w \in \real^{n_{\rm c} \times n_w}$. The goal is to minimize the objective function, subject to the physical limitations imposed by the plant in steady-state, and subject to the engineering constraints imposed by the designer. We make the following standing assumptions on the problem data.

\medskip

\begin{assumption}[\bf Problem Data]\label{Ass:Data}
On the interiors of their domains, the maps $f_0$ and $g_0$ are convex, continuously differentiable, and have at least locally Lipschitz continuous gradients. The set of disturbances $\mathcal{W}$ is convex and compact. The problem \eqref{Eq:OSS} is strictly feasible and has an optimal solution for all $w \in \mathcal{W}$. Both $z(t)$ and the engineering constraint violation $H_z z(t) + H_u u(t) + H_w w(t)$ are measurable.
\end{assumption}

\medskip

Uniqueness of a primal solution to \eqref{Eq:OSS} is guaranteed if $u \mapsto f_0(u) + g_0(G_uu + G_ww)$ is strongly convex on $\mathcal{U} \intersection \setdef{u\in\mathcal{U}}{G_u u + G_w w \in \mathcal{Z}}$ for each $w \in \mathcal{W}$; this holds, for instance, if $f_0$ is strongly convex on $\mathcal{U}$. The dual variable $\mu$ associated with the linear constraint
\[
0 = H_z(G_u\bar{u}+G_ww) + H_u\bar{u} + H_ww
\]
will be unique if $H_zG_u + H_u$ has full row rank; ensuring this is often simply a matter of the designer properly specifying the desired additional constraints \eqref{Eq:OSS-3}.


\subsection{Review of OSS Control}
\label{Sec:PreviousOSS}

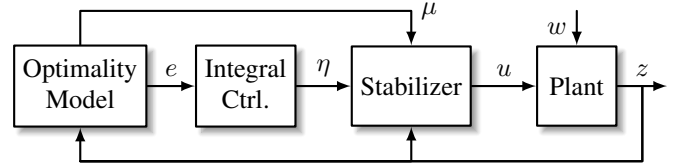
\begin{figure}[ht!]
	\begin{center}  
	    \begin{tikzpicture}[auto, scale = 0.6, node distance=2cm,>=latex', every node/.style={scale=1}]
      \tikzstyle{anch} = [coordinate]
      \tikzstyle{block} = [draw, fill=white, rectangle, 
      minimum height=3em, minimum width=6em, blur shadow={shadow blur steps=5}]
      \tikzstyle{smallblock} = [draw, fill=white, rectangle, 
      minimum height=3em, minimum width=5em, blur shadow={shadow blur steps=5}]
      \tikzstyle{hold} = [draw, fill=white, rectangle, 
      minimum height=3em, minimum width=3em, blur shadow={shadow blur steps=5}]
      \tikzstyle{dzblock} = [draw, fill=white, rectangle, minimum height=3em, minimum width=4em, blur shadow={shadow blur steps=5},
      path picture = {
        \draw[thin, black] ([yshift=-0.1cm]path picture bounding box.north) -- ([yshift=0.1cm]path picture bounding box.south);
        \draw[thin, black] ([xshift=-0.1cm]path picture bounding box.east) -- ([xshift=0.1cm]path picture bounding box.west);
        \draw[very thick, black] ([xshift=-0.5cm]path picture bounding box.east) -- ([xshift=0.5cm]path picture bounding box.west);
        \draw[very thick, black] ([xshift=-0.5cm]path picture bounding box.east) -- ([xshift=-0.1cm, yshift=+0.4cm]path picture bounding box.east);
        \draw[very thick, black] ([xshift=+0.5cm]path picture bounding box.west) -- ([xshift=+0.1cm, yshift=-0.4cm]path picture bounding box.west);
      }]
      \tikzstyle{sum} = [draw, fill=white, circle, node distance=1cm, blur shadow={shadow blur steps=8}]
      \tikzstyle{input} = [coordinate]
      \tikzstyle{output} = [coordinate]
      \tikzstyle{split} = [coordinate]
      \tikzstyle{pinstyle} = [pin edge={to-,thin,black}]
      \node [smallblock] (optmodel) {\makecell[c]{Optimality\\Model}};
      \node [anch, below of=optmodel, node distance=1cm] (anch) {};
      \node [anch, above of=optmodel, node distance=1cm] (anch2) {};
      \node [hold, right of=optmodel, node distance=2.2cm] (integral) {\makecell[c]{Integral\\Ctrl.}};
      \node [hold, right of=integral, node distance=2.2cm] (stabilizer) {Stabilizer};
      \node [hold, right of=stabilizer,
      node distance=2.2cm] (system) {Plant};
      \node [input,name=disturbance, above of = system, node distance=1cm] {};
      \node [output,right of = system, node distance=1.2cm] (output) {};
      \draw [thick, -latex] (optmodel) -- node[name=eps] {$e$} (integral);
      \draw [thick, -latex] (integral) -- node[name=eta] {$\eta$} (stabilizer);
      \draw [thick, -latex] (stabilizer) -- node[name=u] {$u$} (system);
      \draw [thick, -latex] (system) -- node[name=y] {$z$} (output);
      \draw [thick, -latex] (y) |- (anch) -| (optmodel);
      \draw [thick, -latex] (y) |- (anch) -| (stabilizer);
      \draw [thick, -latex] (optmodel.north) |- (anch2) -| node[] {$\mu$} (stabilizer);
      \draw [thick, -latex] (disturbance) -- node[left] {$w$} (system.north);
    \end{tikzpicture}
	\end{center}
	\caption{Block diagram of general OSS control architecture \cite{LSPL-JWSP-EM:18l}.}
	\label{Fig:OSSArch}	
	\end{figure}

In \cite{LSPL-JWSP-EM:18l} a framework was developed for the construction of controllers solving Problem \ref{Prob:OSS}. The key idea introduced was that of an \emph{optimality model}, which is a (potentially, dynamic) nonlinear filter of the form
\begin{equation}\label{Eq:OptimalityModel}
\tau\dot{\mu} = \varphi(\mu,z,u), \qquad e = h(\mu,z,u)
\end{equation}
with state $\mu$, {output $e$}, and time constant $\tau > 0$. The filter \eqref{Eq:OptimalityModel} is said to be an optimality model if the following statement holds: if $(\bar{x},\bar{\mu},\bar{u},\bar{z})$ is an equilibrium of \eqref{Eq:LTI2}, \eqref{Eq:OptimalityModel} satisfying $0 = h(\bar{\mu},\bar{z},\bar{u})$, then $(\bar{z},\bar{u})$ is an optimal solution of \eqref{Eq:OSS}. The idea is that $\dot{\mu}$ and the \emph{error} variable $e$ in \eqref{Eq:OptimalityModel} should together provide a measure of the optimality gap. Specifically, $\dot{\mu} = 0$ along with regulation of $e$ to zero {will} result in in regulation of $(\bar{z},\bar{u})$ to an optimal solution of \eqref{Eq:OSS}, since the plant itself will enforce the constraint \eqref{Eq:OSS-2} in equilibrium. To this end, the error output $e$ from \eqref{Eq:OptimalityModel} is fed to an integral controller
\begin{equation}\label{Eq:IntOSS}
\tau\dot{\eta} = e.
\end{equation}
The optimality model is appended to the output of the plant \eqref{Eq:LTI2}, creating a cascade of the plant, the optimality model, and a bank of integrators. If a stabilizer processing ($z$,$\mu$,$\eta$) can be designed to ensure that the cascade \eqref{Eq:LTI2}, \eqref{Eq:OptimalityModel}, \eqref{Eq:IntOSS} possess a unique exponentially stable equilibrium point, then the optimality model, integral controller, and stabilizer together constitute a solution to the OSS control problem (Figure \ref{Fig:OSSArch}).

In \cite{LSPL-JWSP-EM:18l} several optimality models were designed based on variations of the KKT conditions for \eqref{Eq:OSS}, and results on stabilizability and detectability of the cascaded system \eqref{Eq:LTI2}, \eqref{Eq:OptimalityModel}, \eqref{Eq:IntOSS} were presented for the case of \emph{quadratic} objective functions $f_0, g_0$. When the objective functions are non-quadratic, the design of a stabilizer becomes more challenging, and no guarantees or constructive results were provided in \cite{LSPL-JWSP-EM:18l}. The goal of this paper is to address this gap for the practically-relevant case of an exponentially stable plant.


\section{Development of Optimality Models}
\label{Sec:Optimality}

Under Assumption \ref{Ass:Data} the KKT conditions provide a necessary and sufficient characterization of optimal points of \eqref{Eq:OSS} and can be used to construct optimality models. Different optimality models can be derived by manipulating the optimization problem \eqref{Eq:OSS} and the resulting KKT conditions, and we now pursue this direction.  Section \ref{Sec:Opt1} is conceptually identical to \cite[Prop. 3.3]{LSPL-JWSP-EM:18l}, while Section \ref{Sec:Opt2} contain new results extending the ideas in \cite[Prop. 3.4]{LSPL-JWSP-EM:18l}. Additional results extending the optimality model developed in \cite[Prop. 3.5]{LSPL-JWSP-EM:18l} are omitted due to space limitations, but {are} available in the extended version of this paper \cite{JWSP:22d-extended}.


\subsection{Optimality Model \#1}
\label{Sec:Opt1}

The most obvious initial step is to eliminate the variable $\bar{z}$ from the problem \eqref{Eq:OSS}. Doing so, the resulting Lagrangian is 
\[
\begin{aligned}
L(\bar{u},\bar{\mu}) &= f_0(\bar{u}) + g_0(G_u\bar{u}+G_ww)\\
&\quad + \bar{\mu}^{\T}(H_z(G_u\bar{u}+G_ww) + H_u\bar{u} + H_ww)
\end{aligned}
\]
with multiplier $\bar{\mu} \in \real^{n_{\rm c}}$, leading to the KKT conditions
\begin{equation}\label{Eq:KKT2}
\begin{aligned}
0 &= \nabla f_0(\bar{u}) + G_u^{\T}\nabla g_0(G_u\bar{u}+G_ww)\\
&\qquad + (H_zG_u+H_u)^{\T}\mu\\
0 &= H_{z}(G_u\bar{u}+G_ww) + H_{u}\bar{u} + H_ww.
\end{aligned}
\end{equation}
The approach to constructing an optimality model is now based on replacement of the steady-state output value $\bar{z} = G_u\bar{u} + G_ww$ by the real-time measurement $z(t)$. Following this idea leads to the optimality model
\begin{subequations}\label{Eq:OSOM}
\begin{align}
\label{Eq:OSOM-1}
\tau\dot{\mu} &= H_z{z} + H_u u + H_w w\\
\label{Eq:OSOM-2}
e &= \nabla f_0(u) + G_u^{\T}\nabla g_0(z) + (H_zG_u+H_u)^{\T}\mu
\end{align}
\end{subequations}
where $\tau > 0$ is a tuning parameter. In \cite{LSPL-JWSP-EM:18l}, \eqref{Eq:OSOM} was referred to as the \emph{output subspace} optimality model; integration of the constraint violation in \eqref{Eq:OSOM-1} will ensure primal feasibility, while $e = 0$ in \eqref{Eq:OSOM-2} will ensure {stationarity}.


\subsection{Optimality Model \#2}
\label{Sec:Opt2}

A quite different optimality model can be obtained via a more involved reduction of the linear constraints \eqref{Eq:OSS-2}--\eqref{Eq:OSS-3}. First note that \eqref{Eq:OSS-2}--\eqref{Eq:OSS-3} can be expressed as
\begin{equation}\label{Eq:VectorizeConstraints}
\begin{bmatrix}
I_r & -G_u\\
H_{z} & H_{u}
\end{bmatrix}\begin{bmatrix}
\bar{z} \\ \bar{u}
\end{bmatrix} = \begin{bmatrix}
G_w\\
-H_w
\end{bmatrix}w.
\end{equation}
Let $T = \left[\begin{smallmatrix}T_z \\ T_u\end{smallmatrix}\right] \in \real^{(r+m)\times q}$ be a matrix such that

\begin{equation}\label{Eq:TEquation}
\mathrm{range}(T) = \mathrm{null}\begin{bmatrix}
I_r & -G_u\\
H_z & H_{u}
\end{bmatrix}.
\end{equation}
This allows the constraints  \eqref{Eq:OSS-2}--\eqref{Eq:OSS-3} to be expressed as
\begin{equation}\label{Eq:ConstraintExpanded2}
\begin{bmatrix}
\bar{z} \\ \bar{u}
\end{bmatrix} = \begin{bmatrix}
T_z\\T_u
\end{bmatrix}\bar{\xi} + \begin{bmatrix}z_0(w) \\ u_0(w)\end{bmatrix},
\end{equation}
where $(z_0(w),u_0(w))$ is any solution of \eqref{Eq:VectorizeConstraints} and $\bar{\xi} \in \real^{q}$ is free. The basic idea is that $\mathrm{range}(T)$ captures the subspace on which the physical {steady-state} plant constraints intersect with the specified {steady-state} engineering constraints. 

Using \eqref{Eq:ConstraintExpanded2}, the problem \eqref{Eq:OSS} is now equivalent to the unconstrained minimization problem
\begin{equation}\label{Eq:OSSReduced}
  \begin{aligned}
    \minimize_{\bar{\xi} \in \real^{q}} &\quad f_0(T_u\bar{\xi} + u_0(w)) + g_0(T_z\bar{\xi} + z_0(w))
  \end{aligned}
\end{equation}
with stationarity condition
\begin{equation}\label{Eq:OptimalXi}
0 = T_u^{\T}\nabla f_0(T_u\bar{\xi}+u_0(w)) + T_z^{\T}\nabla g_0(T_z\bar{\xi} + z_0(w)).
\end{equation}
By now re-inserting the real-time signals $u(t)$ and $z(t)$ in place of $\bar{u}$ and $\bar{z}$, we obtain the optimality model
\begin{equation}\label{Eq:RFSError}
e = \begin{bmatrix}
e_1 \\ e_2
\end{bmatrix} = 
\begin{bmatrix}
T_u^{\T}\nabla f_0(u) + T_z^{\T}\nabla g_0(z)\\
H_{z}z + H_{u}u + H_w w
\end{bmatrix}.
\end{equation}
In contrast with \eqref{Eq:OSOM}, \eqref{Eq:RFSError} does not contain a dynamic state corresponding to any dual variable; in \cite{LSPL-JWSP-EM:18l}, a slightly different construction of \eqref{Eq:RFSError} was referred to as the \emph{feasible subspace} optimality model. Regulation of the error signal $e_1 \in \real^{q}$ to zero ensures {stationarity}, while regulation of $e_2 \in \real^{n_{\rm c}}$ to zero ensures primal feasibility. 

If $T$ is selected to have full column rank, then $q = m - n_{\rm c}$, and the total number of error signals in \eqref{Eq:RFSError} is $m - n_{\rm c} + n_{\rm c} = m$, the number of plant inputs. Moreover, in this case, $T_z$ and $T_u$ enjoy several rank properties; the proof of the following can be found in the extended version \cite{JWSP:22d-extended}.

\smallskip

\begin{lemma}\label{Lem:PropT}
If $T = \left[\begin{smallmatrix}T_z \\ T_u\end{smallmatrix}\right] \in \real^{(r+m)\times q}$ is a matrix of full column rank satisfying \eqref{Eq:TEquation}, then the following hold:
\begin{enumerate}
\item \label{Lem:PropT-1} $T_u \in \real^{m \times q}$ has full column rank;
\item \label{Lem:PropT-4} $T_z \in \real^{r \times q}$ has full column rank if and only if 
\begin{equation}\label{Eq:TuCondition}
\mathrm{range}(T_u) \subseteq \mathrm{range}(G_u^{\T}) = \mathrm{range}(G_u^{\dagger});
\end{equation}
\item \label{Lem:PropT-2} $T_z$ has full column rank if $G_u$ has full column rank;
\item \label{Lem:PropT-3} if $G_u$ has full row rank, then $T_z$ may be chosen to have full column rank if and only if there exists a full column rank matrix $X$ such that $(H_{z}+H_{u}G_u^{\dagger})X = 0$, in which case one choice  is $T_z = X$ and $T_u = G_u^{\dagger}X$.
\end{enumerate}
\end{lemma}

%
%
%


\section{Stabilizers for Low-Gain OSS Control}
\label{Sec:ClosedLoop}

This section details stabilizer designs to accompany the optimality models described in Section \ref{Sec:Optimality}. We follow a low-gain integral control approach \cite{EJD:76,JWSP:20a,JWSP:21a}, wherein the controller is tuned to operate slowly compared to the stable plant dynamics \eqref{Eq:LTI2}. Taken together, the optimality model \eqref{Eq:OptimalityModel}, integral controller \eqref{Eq:IntOSS}, and stabilizer (to be designed) will have the nonlinear state-space form
\[
\tau \dot{x}_{\rm c} = f_{\rm c}(x_{\rm c},z,u), \qquad u = k_{\rm c}(x_{\rm c}),
\]
with concatenated state $x_{\rm c}$ and time constant $\tau > 0$. The maps $f_{\rm c}$ and $k_{\rm c}$ will be (at least) locally Lipschitz continuous on the domain of interest. When $\tau$ is large, $x_{\rm c}$ will change slowly, as will $u$, and the plant \eqref{Eq:LTI2} will quickly converge to the quasi-steady-state output $\bar{z}(u,w) \define G_u u + G_ww$. It follows from singular perturbation arguments 
\cite[Chp. 11]{HKK:02} that if the reduced system
\[
\dot{x}_{\rm c} = f_{\rm c}(x_{\rm c},\bar{z}(u,w),u), \qquad u = k_{\rm c}(x_{\rm c}).
\]
possesses a unique locally exponentially stable equilibrium point, then there will exist $\tau^{\star} > 0$ such that for all $\tau \in (0,\tau^{\star})$, the closed-loop system will possess a unique locally exponentially stable equilibrium point, and the controller will solve Problem \ref{Prob:OSS}. In what follows, we therefore jump immediately to analyses of the reduced dynamics arising from our different stabilizer designs, and stability results are stated with the understanding that a full singular perturbation argument can indeed be made.


\subsection{Optimality Model \#1}
\label{Sec:StabOpt1}

For the optimality model \eqref{Eq:OSOM}, we present two approaches for the design of accompanying stabilizers. 

\smallskip


%
%
%

\medskip


\subsubsection{Primal-Dual Stabilizer} 
\label{Sec:PrimalDualStabilizer}

Integrating the error from \eqref{Eq:OSOM} as $\tau\dot{\eta} = -e$, and selecting the simple stabilizer $u = \eta$ leads immediately to the \emph{primal-dual}-type algorithm
\begin{equation}\label{Eq:PrimalDual}
\begin{aligned}
\tau_{\rm p}\dot{u} &= -\nabla f_0(u) - G_u^{\T}\nabla g_0(z) - (H_zG_u+H_u)^{\T}\mu\\
\tau_{\rm d}\dot{\mu} &= H_z{z} + H_u u + H_w w\\
\end{aligned}
\end{equation}
with time constants $\tau_{\rm p}, \tau_{\rm d} > 0$. As the stability of closely related schemes has already been examined in the literature (e.g., \cite{MC-EDA-AB:18}) we simply state the key stability result. 

\smallskip

\begin{proposition}[\bf Primal-Dual]
If $u \mapsto f_0(u) + g_0(G_uu)$ is strongly convex on $\mathcal{U}$, and if $H_zG_u+H_u$ has full row rank, then there exists $\tau^{\star} > 0$ such that for all $(\tau_{\rm p},\tau_{\rm d}) \in (\tau^{\star},\infty)^2$, the controller \eqref{Eq:PrimalDual} solves Problem \ref{Prob:OSS}.
\end{proposition}

\smallskip

\medskip

\subsubsection{Inversion-Based Stabilizer}
\label{Sec:InversionStabilizer}

Our second design departs slightly from the architecture in Figure \ref{Fig:OSSArch} by omitting the bank of integrators. If $\nabla f$ is an invertible mapping, then one may explicitly solve \eqref{Eq:OSOM-2} for $u$ to obtain the controller
\begin{subequations}\label{Eq:Inversion}
\begin{align}
\label{Eq:Inversion-1}
\tau\dot{\mu} &= H_z{z} + H_u u + H_w w\\
\label{Eq:Inversion-2}
u &= (\nabla f_0)^{-1}(-G_u^{\T}\nabla g_0(z) - (H_zG_u+H_u)^{\T}\mu),
\end{align}
\end{subequations}
consisting of an integration on the constraint violation \eqref{Eq:Inversion-1} followed by a nonlinear static output feedback \eqref{Eq:Inversion-2}. 
The following result characterizes {a case} of interest where this is well-posed and leads to closed-loop stability.

\smallskip

\begin{theorem}[\bf Inversion-Based OSS Control]\label{Thm:Inversion}
Consider the LTI system \eqref{Eq:LTI2} with the inversion-based controller \eqref{Eq:Inversion}. Assume that $H_zG_u + H_u$ has full row rank, and that $f_0$ is strongly convex and essentially smooth\footnote{A continuously differentiable convex function $\map{f_0}{\mathcal{U}}{\real}$ is essentially smooth on $\mathcal{U}$ if $f(x) \to \infty$ as $x$ tends to the boundary of $\mathcal{U}$.} on $\mathcal{U}$. 
Then there exists $\tau^{\star} > 0$ such that for all $\tau \in (0,\tau^{\star})$, the controller \eqref{Eq:Inversion} solves Problem \ref{Prob:OSS}.
\end{theorem}

\smallskip

\begin{proof}
The reduced dynamics of the closed-loop \eqref{Eq:LTI2}, \eqref{Eq:Inversion} are obtained by substituting \eqref{Eq:SteadyStateMap} into \eqref{Eq:Inversion}, yielding
\begin{subequations}\label{Eq:ReducedInversion}
\begin{align}\label{Eq:ReducedInversion-1}
\dot{\mu} &= Nu + \tilde{N}w\\
\label{Eq:ReducedInversion-2}
u &= (\nabla f_0)^{-1}(-G_u^{\T}\nabla g_0(G_uu+G_ww) - N^{\T}\mu)
\end{align}
\end{subequations}
where for compactness we have set $N \define H_z G_u+H_u$ and $\tilde{N} \define H_zG_w+H_w$. By the strong convexity and rank assumptions, for each $w \in \mathcal{W}$, \eqref{Eq:OSS} possesses a unique primal-dual optimal solution $(u^{\star},\mu^{\star})$, characterized by \eqref{Eq:KKT2}, and corresponding to the unique equilibrium point of \eqref{Eq:ReducedInversion}. Set $\xi^{\star} \define -G_u^{\T}\nabla g_0(G_uu^{\star}+G_ww) - N^{\T}\mu^{\star}$. Strong convexity and essential smoothness of $f_0$ on $\mathcal{U}$ imply that $\map{(\nabla f_0)^{-1}}{\real^m}{\mathcal{U}}$ is well-defined \cite[Sec. 26]{RTR:70} and globally Lipschitz continuous. Further, since $\nabla f_0$ is locally Lipschitz continuous, it is Lipschitz continuous on a compact set containing $(\nabla f_0)^{-1}(\xi^{\star})$, and it then follows by duality results in \cite{XZ:18} that $(\nabla f_0)^{-1}$ will be strongly monotone and Lipschitz continuous on a compact set containing $\xi^{\star}$. Now define
\[
J_{w}(u) \define f_0(u) + g_0(G_uu+G_ww).
\]
With this, \eqref{Eq:ReducedInversion-2} can be expressed as
\[
\begin{aligned}
\dot{\mu} &= Nu + \tilde{N}w, \qquad \nabla J_w(u) = - N^{\T}\mu.
\end{aligned}
\]
Making the change of variables $\tilde{\mu} = \mu^{\star} - \mu$ and $\tilde{u} = u - u^{\star}$ and using that $(u^{\star},\mu^{\star})$ satisfy \eqref{Eq:KKT2}, we obtain the error dynamics
\[
\dot{\tilde{\mu}} = -N\tilde{u}, \quad \Psi_{w}(\tilde{u}) \define \nabla J_w(\tilde{u}+u^{\star})-\nabla J_w(u^{\star}) = N^{\T}\tilde{\mu}
\]
Since $g_0$ is convex with locally Lipschitz continuous gradient, $\Psi_{w}$ inherits the same properties as $\nabla f_0$. Thus, $\Psi_{w}^{-1}$ is well-defined, and is strongly monotone and locally Lipschtiz on a compact set containing the origin. Eliminating $\tilde{u}$, the reduced dynamics now simplify to 
\[
\dot{\tilde{\mu}} = -N \Psi_{w}^{-1}(N^{\T}\tilde{\mu})
\]
with equilibrium point at $\tilde{\mu} = 0$. With Lyapunov candidate $V(\tilde{\mu}) = \tfrac{1}{2}\|\tilde{\mu}\|_2^2$ we compute that
\[
\begin{aligned}
\dot{V}(\tilde{\mu}) = -(N^{\T}\tilde{\mu})^{\T}\Psi^{-1}_{w}(N^{\T}\tilde{\mu}) &\leq - m_{\Psi} \|N^{\T}\tilde{\mu}\|_2^2
\end{aligned}
\]
for some $m > 0$ and all $\tilde{\mu}$ such that $\|\tilde{\mu}\|_2$ is sufficiently small. Since $N$ has full row rank, we further have $\dot{V}(\tilde{\mu}) \leq -c V(\tilde{\mu})$ for some $c > 0$, which establishes local exponential stability of the reduced dynamics.
\end{proof}


\subsection{Optimality Model \#2}
\label{Sec:StabOpt2}

For the optimality model \eqref{Eq:RFSError}, we present two approaches for the design of accompanying stabilizers. 

\smallskip

\subsubsection{Explicit Two-Loop Stabilizer Design}
\label{Sec:TwoLoop}

\begin{figure}[ht!]
	\begin{center}  
	    \begin{tikzpicture}[auto, scale = 0.6, node distance=2cm,>=latex', every node/.style={scale=1}]
      \tikzstyle{anch} = [coordinate]
      \tikzstyle{block} = [draw, fill=white, rectangle, 
      minimum height=3em, minimum width=6em, blur shadow={shadow blur steps=5}]
      \tikzstyle{smallblock} = [draw, fill=white, rectangle, 
      minimum height=3em, minimum width=5em, blur shadow={shadow blur steps=5}]
      \tikzstyle{thinblock} = [draw, fill=white, rectangle, 
      minimum height=2em, minimum width=5em, blur shadow={shadow blur steps=5}]
      \tikzstyle{hold} = [draw, fill=white, rectangle, 
      minimum height=3em, minimum width=3em, blur shadow={shadow blur steps=5}]
      \tikzstyle{smallhold} = [draw, fill=white, rectangle, 
      minimum height=2em, minimum width=2em, blur shadow={shadow blur steps=5}]
      \tikzstyle{vsmallhold} = [draw, fill=white, rectangle, 
      minimum height=1.5em, minimum width=1.5em, blur shadow={shadow blur steps=5}]
      \tikzstyle{dzblock} = [draw, fill=white, rectangle, minimum height=3em, minimum width=4em, blur shadow={shadow blur steps=5},
      path picture = {
        \draw[thin, black] ([yshift=-0.1cm]path picture bounding box.north) -- ([yshift=0.1cm]path picture bounding box.south);
        \draw[thin, black] ([xshift=-0.1cm]path picture bounding box.east) -- ([xshift=0.1cm]path picture bounding box.west);
        \draw[very thick, black] ([xshift=-0.5cm]path picture bounding box.east) -- ([xshift=0.5cm]path picture bounding box.west);
        \draw[very thick, black] ([xshift=-0.5cm]path picture bounding box.east) -- ([xshift=-0.1cm, yshift=+0.4cm]path picture bounding box.east);
        \draw[very thick, black] ([xshift=+0.5cm]path picture bounding box.west) -- ([xshift=+0.1cm, yshift=-0.4cm]path picture bounding box.west);
      }]
      \tikzstyle{sum} = [draw, fill=white, circle, node distance=1cm, blur shadow={shadow blur steps=8}]
      \tikzstyle{input} = [coordinate]
      \tikzstyle{output} = [coordinate]
      \tikzstyle{split} = [coordinate]
      \tikzstyle{pinstyle} = [pin edge={to-,thin,black}]
      \node [smallhold] (system) {Plant};
      \node [smallblock, right of=system, node distance=2.2cm] (optmodel) {\makecell[c]{Optimality\\Model \eqref{Eq:RFSError}}};
      \draw [thick, -latex] (system) -- node[] {$z$} (optmodel);
      \node [sum, left of=system, node distance=1.5cm] (sum2) {};
      \node [vsmallhold, below of=sum2, node distance=0.8cm] (K2) {$K_2$};
      
      \node [sum, left of=system, node distance=2.5cm] (sum1) {};
      \node [vsmallhold, below of=sum1, node distance=0.8cm] (K1) {$K_1$};
      
      \node [thinblock, below of=system, node distance=1.4cm, xshift=0.4cm] (int1) {$\tau_2\dot{\eta}_2 = -e_2$};
      \draw [thick, -latex] (int1) -| (K2) -- (sum2);
      
      \node [thinblock, below of=int1, node distance=1cm] (int2) {$\tau_1\dot{\eta}_1 = -e_1$};
      \draw [thick, -latex] (int2) -| (K1) -- (sum1);
      
      \draw [thick, -latex] (sum1) -- node[] {$u_1$} (sum2);
      \draw [thick, -latex] (sum2) -- node[] {$u$} (system);
      
      \draw [thick, -latex] (optmodel) |- (int1);
       \draw [thick, -latex] (optmodel) |- (int2);
      
      \node [input, name=disturbance, above of = system, node distance=1cm] {};
      \draw [thick, -latex] (disturbance) -- node[left] {$w$} (system.north);
%
    \end{tikzpicture}
	\end{center}
	\caption{Two-loop stabilizer design for OM \#2.}
	\label{Fig:TwoLoop}	
	\end{figure}
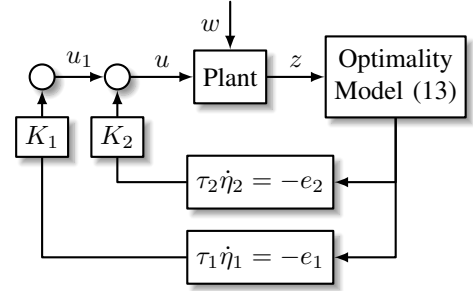

Integrating both error signals from the optimality model \eqref{Eq:RFSError} with two time constants $\tau_1,\tau_2 > 0$, we look for an integral feedback design of the form
\begin{equation}\label{Eq:RFS-Stab}
\begin{aligned}
\tau_1 \dot{\eta}_1 &= -e_1 = -T_u^{\T}\nabla f_0(u) - T_z^{\T}\nabla g_0(z)\\
\tau_2 \dot{\eta}_2 &= -e_2 = -(H_{z}z + H_{u}u + H_w w)\\
u &= K_1\eta_1 + K_2\eta_2
\end{aligned}
\end{equation}
for matrices $K_1, K_2$ to be designed. Inspired by \cite{EJD:76}, the idea we will pursue is the sequential closing of these two integral control loops. First, the loop involving $\eta_2$ will be closed and tuned assuming $\eta_1$ is absent, then the $\eta_1$ loop will be closed around the $\eta_2$ loop, as show in Figure \ref{Fig:TwoLoop}.

The following result provides tuning criteria.

\smallskip

\begin{theorem}[\bf Two-Loop Stabilizer]\label{Thm:TwoLoop}
Consider the plant \eqref{Eq:LTI2} with controller \eqref{Eq:RFS-Stab}, and assume that $T = \left[\begin{smallmatrix}T_z \\ T_u\end{smallmatrix}\right]$ has been selected to have full column rank. 
Suppose that $N \define H_zG_u + H_u$ has full row rank and that $\xi \mapsto f_0(T_u\xi + u_0(w)) + g_0(T_z\xi + z_0(w))$ is strongly convex on $\Xi_{w} = \setdef{\xi}{T_u\xi + u_0(w) \in \mathcal{U}}$ for each $w \in \mathcal{W}$. Select $K_2$ such that $-NK_2$ is Hurwitz, define the projection matrix
\begin{equation}\label{Eq:Pic}
\Pi_{\rm c} \define I_m-K_2[NK_2]^{-1}N,
\end{equation}
{and finally, for any $P \succ 0$, select} $K_1$ of full column rank such that $\Pi_{\rm c}K_1 = T_uP$. Then the controller \eqref{Eq:RFS-Stab} solves Problem \ref{Prob:OSS} for all tunings of the form $\tau_1 \gg \tau_2 \gg 0$. 
\end{theorem}

The convexity assumptions required in Theorem \ref{Thm:TwoLoop} are slightly weaker than those in the results of Section \ref{Sec:StabOpt1}. Specifically, the present convexity condition is a statement about the behaviour of the objective functions on the affine solution set to \eqref{Eq:VectorizeConstraints}, as opposed to, e.g., $f_0$ being strongly convex on its entire domain $\mathcal{U}$.

\begin{pfof}{Theorem \ref{Thm:TwoLoop}}
The proof is based on two sequential applications of \cite[Thm 3.1]{JWSP:20a}; the constructions for the controller gains $K_1, K_2$ will be shown to be well-posed along the way. First, consider the controller
\begin{equation}\label{Eq:FirstController}
\tau_2\dot{\eta}_2 = -(H_{z}z + H_{u}u + H_w w), \quad u = K_2\eta_2 + u_1
\end{equation}
with \emph{constant} auxiliary input $u_1$. With the plant \eqref{Eq:LTI2} interconnected with \eqref{Eq:FirstController}, the reduced dynamics will have the form
\begin{equation}\label{Eq:Red1}
\dot{\eta}_2 = -Nu - (H_zG_w+H_w) w, \quad u = K_2\eta_2 + u_1.
\end{equation}
Since $N$ has full row rank, there exists $K_2$ such that $-NK_2$ is Hurwitz; for instance, take $K_2 = N^{\dagger} = N^{\T}(NN^{\T})^{-1}$. It follows that the reduced dynamics \eqref{Eq:Red1} are globally exponentially stable, and we conclude from \cite[Thm 3.1]{JWSP:20a} that we may select $\tau_2 > 0$ sufficiently large such that the (LTI) closed-loop system \eqref{Eq:LTI2} with controller \eqref{Eq:FirstController} is internally exponentially stable and achieves $e_2(t) \to 0$ as $t \to \infty$ for all asymptotically constant exogenous inputs $w$ and $u_1$. 

We now consider the plant defined by \eqref{Eq:LTI2} and \eqref{Eq:FirstController} with state $(x,\eta_2)$, inputs $(u_1,w)$ and nonlinear output $e_1$ defined in \eqref{Eq:RFSError}. By the previous constructions, the state dynamics are LTI and internally exponentially stable. We {seek} to compute the steady-state input-output relation of this system, which is characterized by the equations
\[
\begin{aligned}
\bar{z} &= G_u \bar{u} + G_w w, \,\, &0 &= N\bar{u} + (H_zG_w + H_w)w\\
\bar{u} &= K_2\bar{\eta}_2 + \bar{u}_1,\,\, &\bar{e}_1 &= T_u^{\T}\nabla f_0(\bar{u}) + T_z^{\T}\nabla g_0(\bar{z}).
\end{aligned}
\]
Solving, one finds that
\[
\bar{\eta}_2 = -[NK_2]^{-1}N\bar{u}_1 - \mathcal{D}w
\]
where $\mathcal{D} = [NK_2]^{-1}(H_zG_w + H_w)$, and therefore that
\[
\begin{aligned}
\bar{u} &= \bar{u}_1 - K_2[NK_2]^{-1}N\bar{u}_1 - K_2\mathcal{D}w\\
&= \Pi_{\rm c}\bar{u}_1 - K_2\mathcal{D}w.
\end{aligned}
\]
where $\Pi_{\rm c}$ is as in \eqref{Eq:Pic}. A straightforward calculation shows that $\Pi_{\rm c}^2 = \Pi_{\rm c}$, so $\Pi_{\rm c}$ is an oblique projection matrix. Moreover, for later use, note that $\mathrm{null}(N) \subseteq \mathrm{range}(\Pi_{\rm c})$. Substituting our expression for $\bar{u}$ back into the expression for $\bar{e}_1$, we obtain the steady-state input-output relationship
\begin{equation}\label{Eq:pi2}
\begin{aligned}
\bar{e}_1 &= T_u^{\T}\nabla f_0(\Pi_{\rm c}\bar{u}_1 + u_0(w))\\
&\qquad \qquad + T_z^{\T}\nabla g_0(G_u\Pi_{\rm c}\bar{u}_1 + z_0(w))
\end{aligned}
\end{equation}
where $(u_0(w),z_0(w))$ are as in \eqref{Eq:ConstraintExpanded2}. We now consider the feedback controller
\begin{equation}\label{Eq:SecondController}
\tau_1\dot{\eta}_1 = -e_1, \qquad u_1 = K_1\eta_1,
\end{equation}
where $K_1$ is selected such that $\Pi_{\rm c}K_1 = T_u P$ with $P \succ 0$. From Lemma \ref{Lem:PropT} \ref{Lem:PropT-1} we know that $T_u$ has full column rank. Moreover, by construction from \eqref{Eq:TEquation}, we have that $(H_z G_u + H_u)T_u = 0$, meaning that $\mathrm{range}(T_u) \subseteq \mathrm{null}(N) \subseteq \mathrm{range}(\Pi_{\rm c})$. Thus, the equation $\Pi_{\rm c}K_1 = T_u P$ is always solvable for a matrix $K_1$ of full column rank. To complete the proof, we verify that the reduced dynamics, given by
\[
\dot{\eta}_1 = -T_u^{\T}\nabla f_0(\Pi_{\rm c}K_1\eta_1 + u_0) - T_z^{\T}\nabla g_0(G_u\Pi_{\rm c}K_1\eta_1 + z_0)
\]
are exponentially stable. Substituting for $\Pi_{\rm c}K_1$ and using from \eqref{Eq:TEquation} that $T_z = G_uT_u$, we can simplify this to obtain
\[
\begin{aligned}
\dot{\eta}_1 &= -T_u^{\T}\nabla f_0(T_u P\eta_1 + u_0) - T_z^{\T}\nabla g_0(G_u T_u P\eta_1 + z_0)\\
&= -T_u^{\T}\nabla f_0(T_u P\eta_1 + u_0) - T_z^{\T}\nabla g_0(T_z P\eta_1 + z_0)\\
&\define -\Phi(P\eta_1)
\end{aligned}
\]
The main convexity assumption implies that $\Phi$ is strongly monotone on the feasible set. Similar to the proof of Theorem \ref{Thm:Inversion}, a simple Lyapunov analysis with candidate $V(\eta_1) = \|\eta_1-\eta_1^{\star}\|_{P}^2$ now completes the proof of local exponential stability of the reduced dynamics.
\end{pfof}


\subsubsection{Robust \& Optimal Stabilizer Design}
\label{Sec:RobustDesign}

Theorem \ref{Thm:TwoLoop} provides an explicit design for an appropriate gain $K = \left[\begin{smallmatrix}K_1 & K_2\end{smallmatrix}\right]$ for use in \eqref{Eq:RFS-Stab}. The quality of this design can be significantly improved by leveraging techniques from robust and optimal control theory. With the plant in quasi steady-state as given by \eqref{Eq:SteadyStateMap}, the reduced dynamics associated with \eqref{Eq:RFS-Stab} will have the form
\begin{subequations}\label{Eq:CLSRobust}
\begin{align}
\label{Eq:CLSRobust-1}
\dot{\eta}_1 &= -T_u^{\T}\nabla f_0(u) - T_z^{\T}\nabla g_0(G_uu + G_w w)\\
\label{Eq:CLSRobust-2}
\dot{\eta}_2 &= -Nu - \tilde{N}w\\
\label{Eq:CLSRobust-3}
u &= K_1\eta_1 + K_2\eta_2.
\end{align}
\end{subequations}
where, as before, we set $N \define H_zG_u + H_u$ and $\tilde{N} \define H_zG_w+H_w$. We now recognize the design of $K_1, K_2$ in \eqref{Eq:CLSRobust} a \emph{state-feedback} design problem. 
To the system \eqref{Eq:CLSRobust}, we associate the performance outputs $z_1 = \dot{\eta}_1$ and $z_2 = \rho \dot{\eta}_2$ corresponding to the stationarity and primal feasibility violations, respectively, with $\rho > 0$ a tuning parameter. To minimize conservatism in the design that follows, we will separate out any quadratic cost portions from $f_0$ and $g_0$, so that
\[
\begin{aligned}
\nabla f_0(u) &= Q_1 u + \nabla \tilde{f}_0(u), \quad \nabla g_0(z) = Q_2 z + \nabla \tilde{g}_0(z)
\end{aligned}
\]
for matrices $Q_1, Q_2 \succeq 0$ and $\tilde{f_0},\tilde{g}_0$ being convex. With this, the overall system can be represented in standard linear-fractional representation (LFR) form \cite{LEG-SN:00,CS-SW:15} as
%
\[
\def\arraystretch{1.2}
\begin{aligned}
\left[\begin{array}{@{}c@{}}
\dot{\eta}_1\\
\dot{\eta}_2\\
\hline
q_1\\
q_2\\
\hline
z_1\\
z_2
\end{array}\right] &= 
\left[\begin{array}{@{}cc|c|cc|c@{}}
0 & 0 & -Q & -T_u^{\T} & -T_z^{\T} & -T_z^{\T}Q_2G_w\\
0 & 0 & -N & 0 & 0 & -\tilde{N}\\
\hline
0 & 0 & I & 0 & 0 & 0\\
0 & 0 & G_u & 0 & 0 & G_w\\
\hline
0 & 0 & -Q & -T_u^{\T} & -T_z^{\T} & -T_z^{\T}Q_2G_w\\
0 & 0 & -\rho N & 0 & 0 & -\rho\tilde{N}\\
\end{array}\right]\left[\begin{array}{@{}c@{}}
\eta_1 \\ \eta_2\\ \hline u\\ \hline p_1 \\ p_2\\ \hline w
\end{array}\right]\\
&\define \left[\begin{array}{@{}c|c|c|c@{}}
\mathcal{A} & \mathcal{B} & \mathcal{B}_1 & \mathcal{B}_2\\
\hline
\mathcal{C}_1 & \mathcal{E}_1 & \mathcal{D}_1 & \mathcal{D}_{12}\\
\hline
\mathcal{C}_2 & \mathcal{E}_2 & \mathcal{D}_{21} & \mathcal{D}_{2}\\
\end{array}\right]\mathrm{col}(\eta_1,\eta_2,u,p_1,p_2,w)
\end{aligned}
\]
where $Q \define T_u^{\T}Q_1 + T_{z}^{\T}Q_2G_u$, and with
\[
\begin{aligned}
p_1 &= \Delta_1(q_1) \define \nabla \tilde{f}_0(q_1)\\
p_2 &= \Delta_2(q_2) \define \nabla \tilde{g}_0(q_2).
\end{aligned}
\]
Since $\nabla f_0$ and $\nabla g_0$ are gradients of convex functions, the blocks $\Delta_i$ satisfy (incremental) pointwise sector constraints of the form \cite{LL-BR-AP:16}
\begin{equation}\label{Eq:IQC}
\begin{bmatrix}
\Delta_i(q_i) - \Delta_i(q_i^{\prime})\\
q_i - q_i^{\prime}
\end{bmatrix}^{\T}\Theta_i\begin{bmatrix}
\Delta(q_i) - \Delta(q_i^{\prime})\\
q_i - q_i^{\prime}
\end{bmatrix} \geq 0
\end{equation}
for all arguments $q_i,q_i^{\prime}$ and for some symmetric block two-by-two matrix $\Theta_i$. The specific form of $\Theta_i$ depends on whether the objective function is (i) convex, (ii) convex with Lipschitz continuous gradient, (iii) strongly convex, or (iv) strongly convex and with Lipschitz continuous gradient. For example, if $f_0$ is strongly convex with parameter $m_f$ and its gradient has Lipschitz constant $L_f$, we may take \cite{LL-BR-AP:16}
\[
\Theta_1 = \begin{bmatrix}
-2 & m_f + L_f\\
m_f + L_f & -2m_f L_f
\end{bmatrix} \otimes I_m.
\]
With this in mind, for selected $\Theta_1, \Theta_2$, we set
\[
\begin{aligned}
\boldsymbol{\Theta} &\define \setdef{\mathrm{daug}(\Theta_1/\theta_1,\Theta_2/\theta_2)}{\theta_1,\theta_2 > 0}\\
\Theta_{\gamma} &\define \mathrm{diag}(-I_{n_w},\tfrac{1}{\gamma^2}I_{q+n_{\rm c}}).
\end{aligned}
\]
with $\gamma > 0$, and where $\mathrm{daug}$ denotes the diagonal augmentation operation \cite{LEG-SN:00}. By standard arguments involving the $S$-Procedure, the state feedback $u = K\left[\begin{smallmatrix}\eta_1 \\ \eta_2\end{smallmatrix}\right]$ will render the reduced dynamics \eqref{Eq:CLSRobust} exponentially stable and achieve $\mathscr{L}_2$-performance on the $w \mapsto z$ channel strictly upper bounded by $\gamma$ if there exists $P \succ 0$ and $\Theta \in \boldsymbol{\Theta}$ such that
\begin{equation}\label{Eq:RobustPerf}
{\small
\def\arraystretch{1.2}
\begin{aligned}
&(\star)^{\sf T}
\left[\begin{array}{@{}cc|c|c@{}}
0 & P & 0 & 0\\
P & 0 & 0 & 0\\
\hline
0 & 0 & \Theta & 0\\
\hline
0 & 0 & 0 & \Theta_{\gamma}
\end{array}\right]
\left[\begin{array}{@{}ccc@{}}
I_{q+n_{\rm c}} & 0 & 0\\
\mathcal{A}+\mathcal{B}K & \mathcal{B}_1 & \mathcal{B}_2\\
\hline
0 & I_{m+r} & 0\\
\mathcal{C}_1+\mathcal{E}_1K & \mathcal{D}_1 & \mathcal{D}_{12}\\
\hline
0 & 0 & I_{n_w}\\
\mathcal{C}_2+\mathcal{E}_2K & \mathcal{D}_{21} & \mathcal{D}_2\\
\end{array}\right] \prec 0.
\end{aligned}
}
\end{equation}
Defining $Y = P^{-1} \succ 0$, performing a congruence transformation on \eqref{Eq:RobustPerf} using $\mathrm{diag}(Y,I,I)$, and applying the \emph{Dualization Lemma} \cite[Cor. 4.10]{CS-SW:15}, \eqref{Eq:RobustPerf} is equivalent to the LMI problem \eqref{Eq:RobSynDual} of finding $Y \succ 0$, $Z \in \real^{m \times (q+n_{\rm c})}$, and $\theta_1,\theta_2 > 0$ such that
\begin{figure*}[t!]
\begin{equation}\label{Eq:RobSynDual}
{\small
(\star)^{\sf T}
\left[\begin{array}{@{}cc|c|cc@{}}
0 & -I_{q+n_{\rm c}} & 0 & 0 & 0\\
-I_{q+n_{\rm c}} & 0 & 0 & 0 & 0\\
\hline
0 & 0 & -\Theta^{-1} & 0 & 0\\
\hline
0 & 0 & 0  & I_{n_w} & 0\\
0 & 0 & 0 & 0 & -\gamma^2 I_{m + n_{\rm c}}\\
\end{array}\right]
\def\arraystretch{1.15}
\left[\begin{array}{@{}ccc@{}}
-(\mathcal{A}Y+\mathcal{B}Z)^{\T} & -(\mathcal{C}_1Y+\mathcal{E}_1Z)^{\T} & -(\mathcal{C}_2Y+\mathcal{E}_2Z)^{\T}\\
I_{q+n_{\rm c}} & 0 & 0\\
\hline
-\mathcal{B}_1^{\T} & -\mathcal{D}_1^{\T} & -\mathcal{D}_{21}^{\T}\\
0 & I_{m+r} & 0\\
\hline
-\mathcal{B}_2^{\T} & -\mathcal{D}_{12}^{\T} & -\mathcal{D}_{2}^{\T}\\
0 & 0 & I_{q+n_{\rm c}}
\end{array}\right] \prec 0.
}
\end{equation}
\end{figure*}
The resulting feedback gain is then recovered as $K = ZY^{-1}$. A design minimizing $\gamma$ can then be obtained via semidefinite programming by minimizing $\gamma^2$ subject to \eqref{Eq:RobSynDual}. We summarize in the following result.

\smallskip

\begin{proposition}[\bf Robust Optimal Stabilizer]\label{Prop:Optimal}
Suppose that \eqref{Eq:RobSynDual} is feasible, and set $K = ZY^{-1}$. Then the controller \eqref{Eq:RFS-Stab} solves Problem \ref{Prob:OSS} for all sufficiently large tunings of the form $\tau_1 = \tau_2 = \tau \gg 0$. 
\end{proposition}

\smallskip

The interested reader will have no issues extending the robust/optimal design approach above to other optimality models; the details are omitted. The only pieces of plant model information required for this LMI-based design procedure are the plant DC gain matrices $G_u$ and $G_w$.



\section{Examples}
\label{Sec:Simulation}

\subsection{Simulation Example}

We first illustrate our ideas via an academic example. The LTI system \eqref{Eq:LTI2} is a randomly generated underactuated stable system with $x \in \real^{30}$, $u \in \real^4$, and $z \in \real^5$, and is such that $G_u$ has full column rank. The problem of interest is asymptotic tracking of step reference signals $r_1,r_2$ for $z_1, z_2$, while constraining $z_3,z_4,z_5$ within specified limits. This should be achieved with minimum control effort, and subject to hard constraints on controls $u$. We formulate this as
\begin{subequations}
  \begin{align}
  \label{Eq:Servo-1}
    \minimize &\quad  \left[\sum_{k=1}^{m}\nolimits \tfrac{1}{2}\bar{u}_k^2 + \gamma \mathsf{B}(\bar{u}_k)\right] + c\,\mathsf{P}(\bar{z})\\
    \label{Eq:Servo-2}
    \subto &\quad \bar{z} = G_u\bar{u} + G_ww\\
    \label{Eq:Servo-3}
           &\quad 0 = z_i - r_i, \qquad i \in \{1,2\}
  \end{align}
\end{subequations}
where $\gamma, c > 0$ are design parameters for the barrier and penalty functions
\[
\begin{aligned}
\mathsf{B}(u_k) &= -\log(u_k^{\rm max}-u_k) - \log(-u_{k}^{\rm min}+u_k)\\
\mathsf{P}(\bar{z}) &= \tfrac{1}{2}\sum_{k=3}^{5}\nolimits \max(0,z_k^{\rm min}-\bar{z}_k,\bar{z}_k-z_{k}^{\rm max})^2
\end{aligned}
\]
Note that the objective function $f_0$ is strongly convex and essentially smooth on $\mathcal{U} = \prod_{k} [u_{k}^{\rm min},u_{k}^{\rm max}]$.

\medskip



\begin{figure}[ht!]
\centering
\begin{subfigure}{0.99\linewidth}
\includegraphics[width=\linewidth]{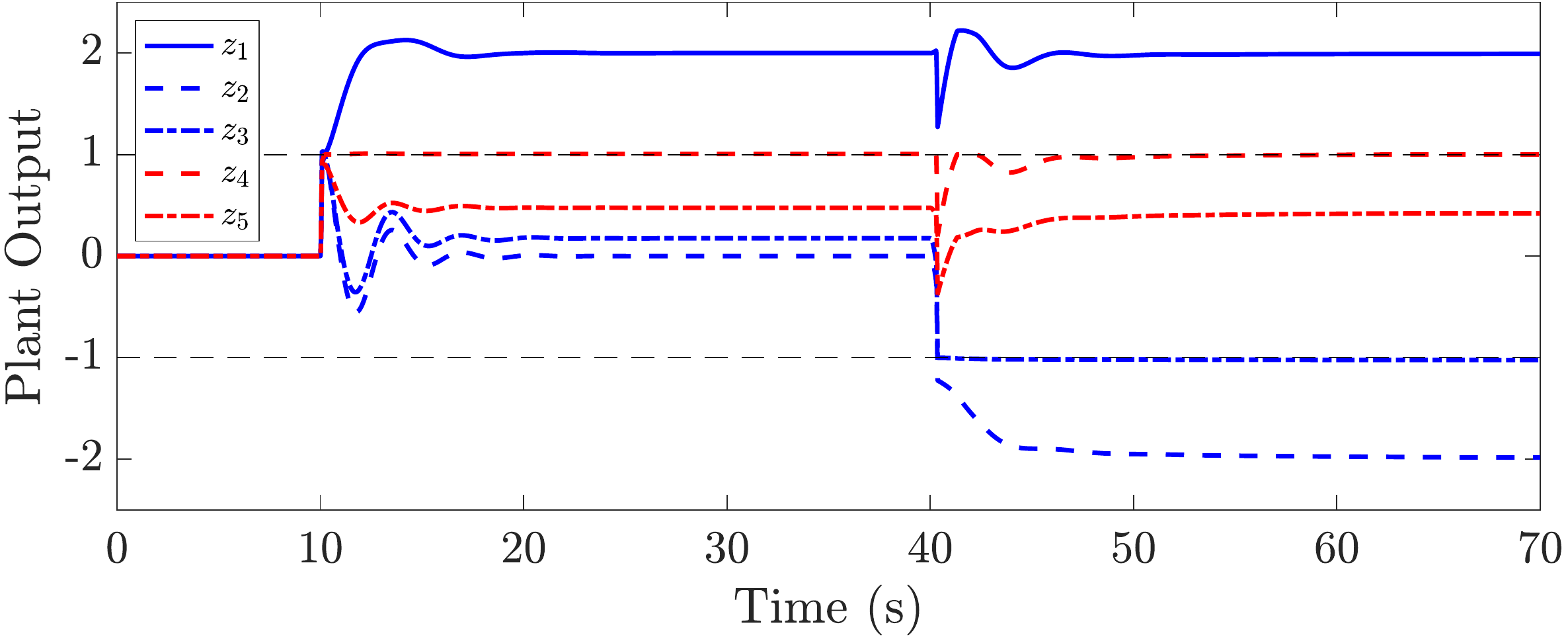}
\label{Fig:OSS-Inv-z}
\end{subfigure}
\begin{subfigure}{0.99\linewidth}
\includegraphics[width=\linewidth]{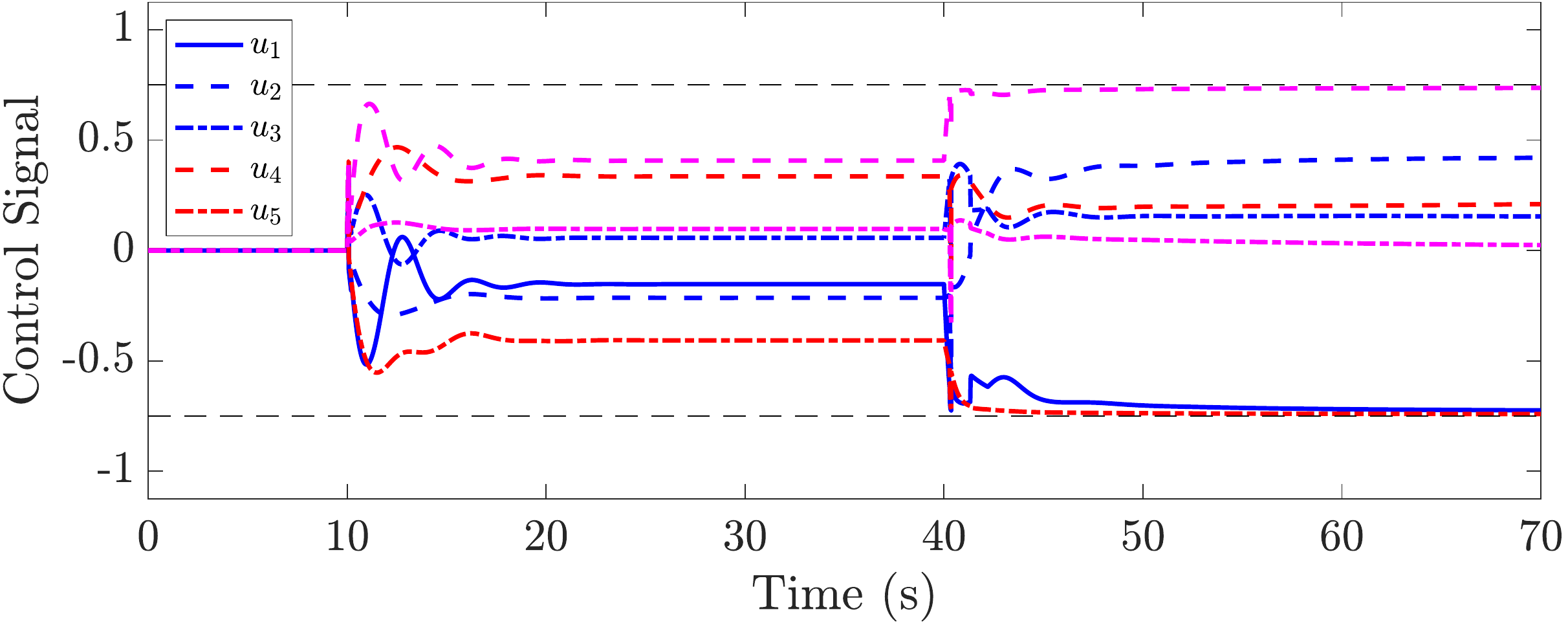}
\label{Fig:OSS-Inv-u}
\end{subfigure}
\caption{System with inversion-based stabilizer; $\tau = 2$.}
\label{Fig:OSS-Inv}
\end{figure}


\begin{figure}[ht!]
\centering
\begin{subfigure}{0.99\linewidth}
\includegraphics[width=\linewidth]{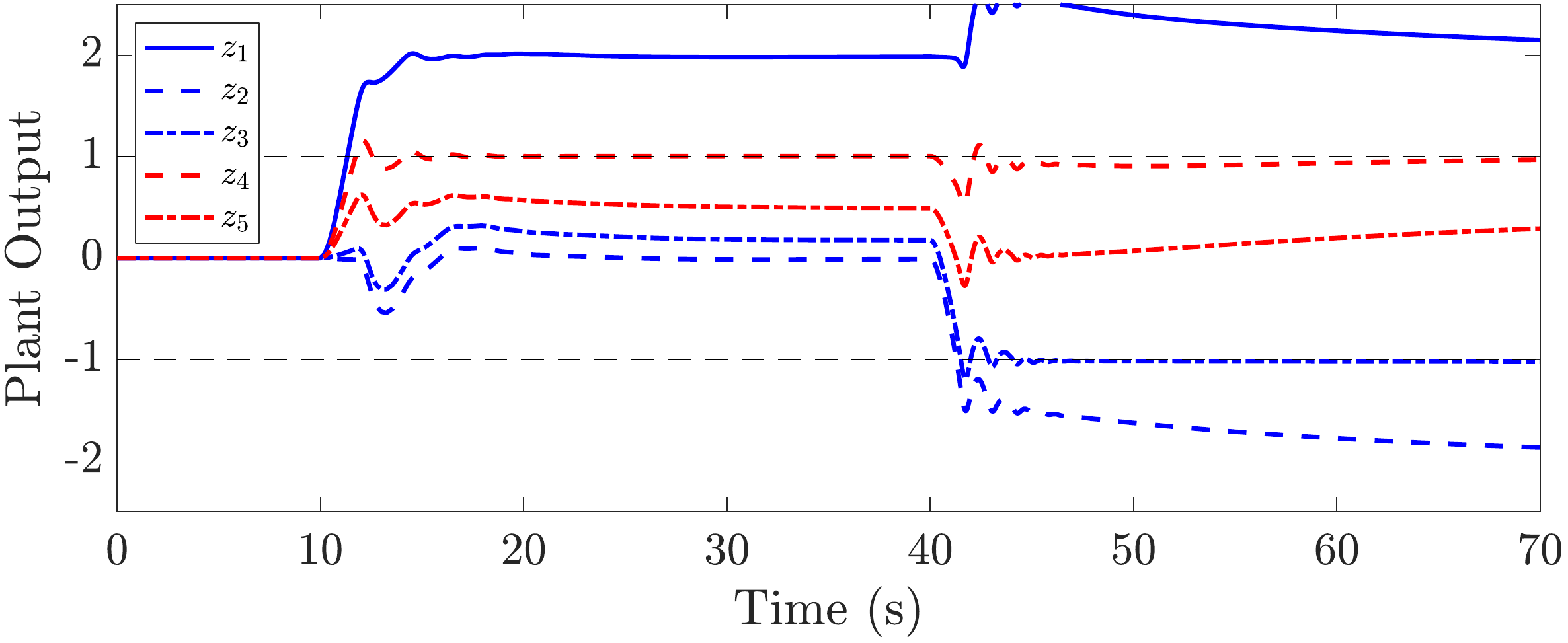}
\label{Fig:OSS-TwoLoop-z}
\end{subfigure}
\begin{subfigure}{0.99\linewidth}
\includegraphics[width=\linewidth]{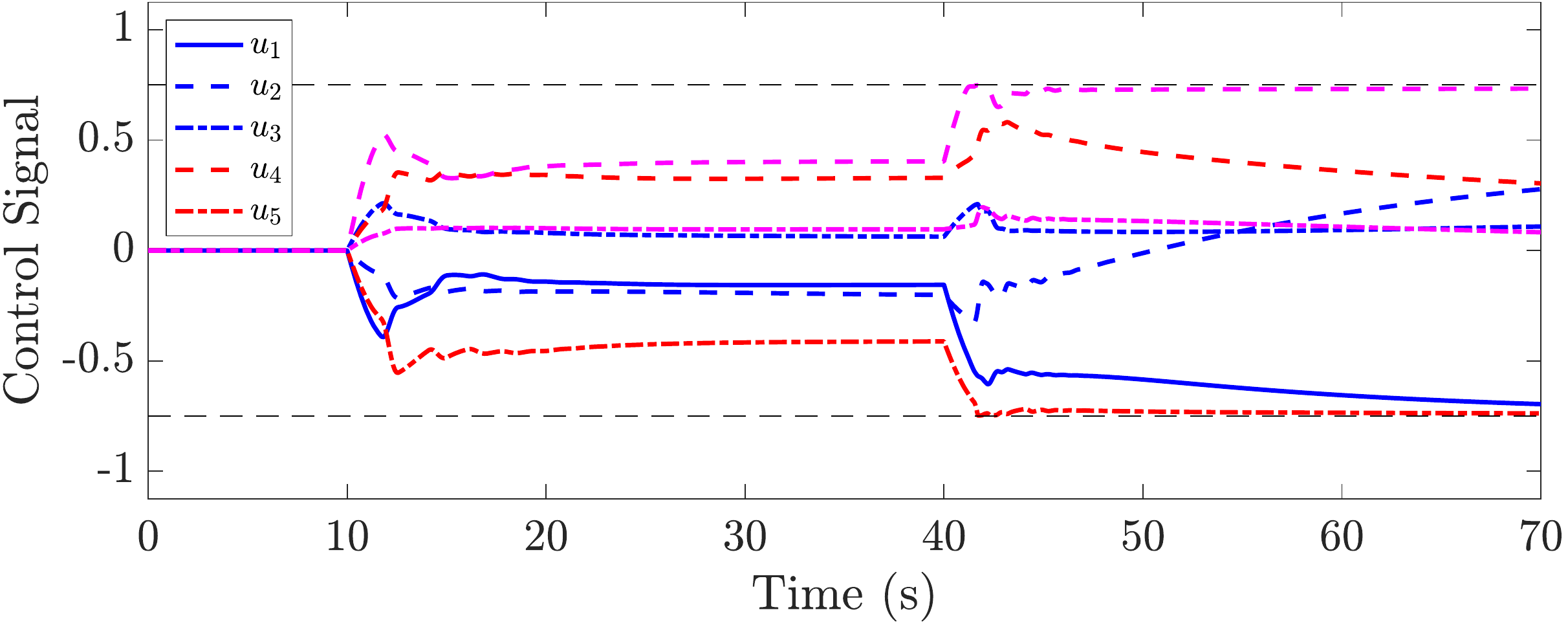}
\label{Fig:OSS-TwoLoop-u}
\end{subfigure}
\caption{System with two-loop stabilizer; $\tau_1 = 5, \tau_2 = 1$.}
\label{Fig:OSS-TwoLoop}
\end{figure}

\begin{figure}[ht!]
\centering
\begin{subfigure}{0.99\linewidth}
\includegraphics[width=\linewidth]{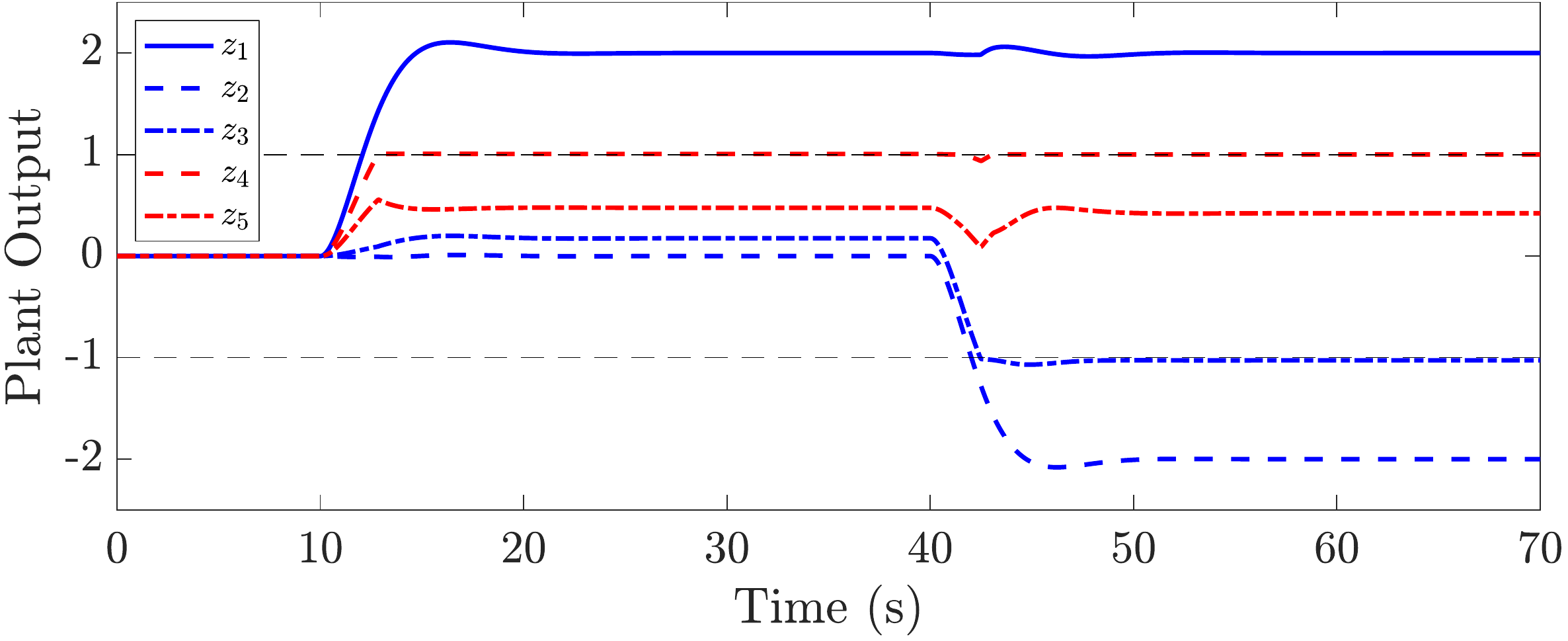}
\label{Fig:OSS-TwoLoop-Optimal-z}
\end{subfigure}
\begin{subfigure}{0.99\linewidth}
\includegraphics[width=\linewidth]{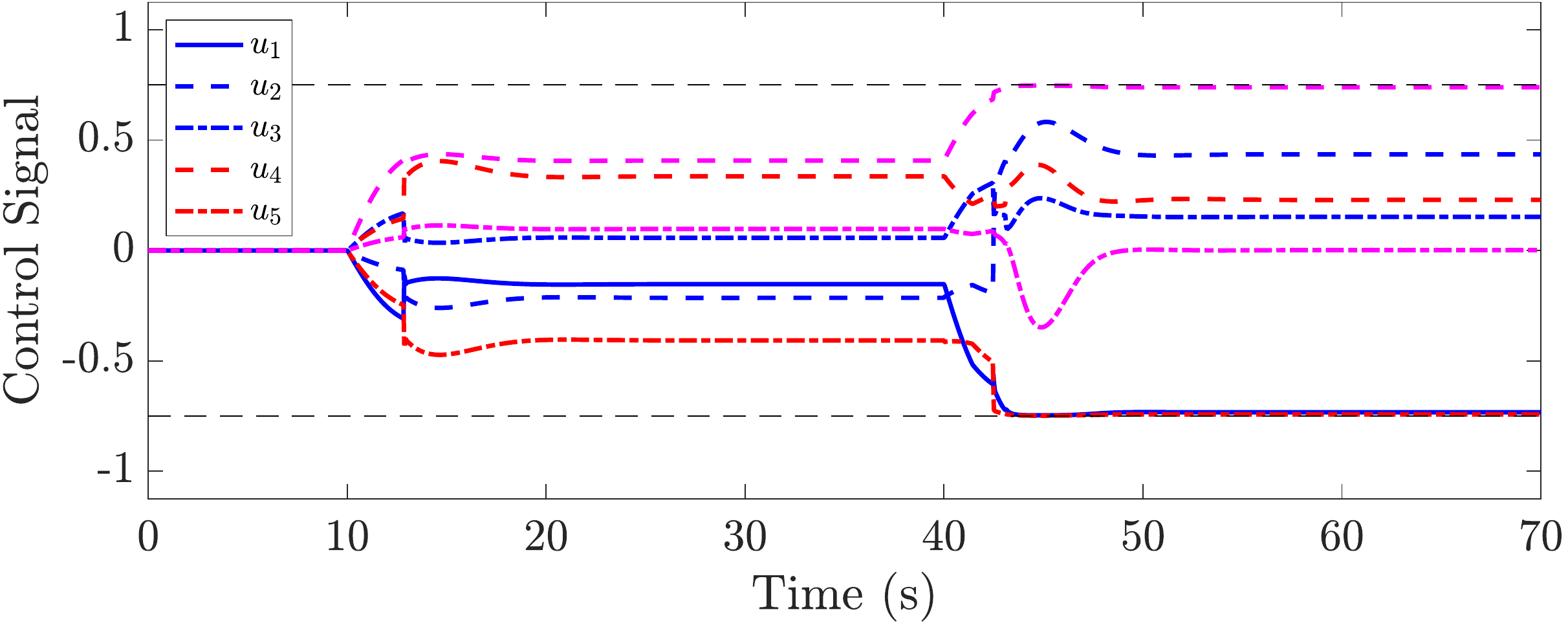}
\label{Fig:OSS-TwoLoop-Optimal-u}
\end{subfigure}
\caption{System with optimal stabilizer; $\rho = 100, \tau = 0.02$.}
\label{Fig:OSS-TwoLoop-Optimal}
\end{figure}

In the following tests we set $u_k^{\rm min} = -0.75$, $u_k^{\rm max} = 0.75$, $z_k^{\rm min} = -1$, $z_k^{\rm max} = 1$ for all $k$, $\gamma = 0.01$, and $c = 50$. The response of the closed-loop system to sequential step changes in the references $(r_1,r_2) = (2,-2)$ at $t = 10$s and $t=40$s is simulated. We illustrate the performance of the inversion-based controller \eqref{Eq:Inversion} in Figure \ref{Fig:OSS-Inv}, the two-loop controller \eqref{Eq:RFS-Stab} in Figure \ref{Fig:OSS-TwoLoop}, and the optimal controller of Proposition \ref{Prop:Optimal} in Figure \ref{Fig:OSS-TwoLoop-Optimal}.\footnote{The primal-dual controller \eqref{Eq:PrimalDual} produces results similar to the inversion-based controller when the time constant $\tau_{\rm p}$ is small compared to $\tau_{\rm d}$.} In all cases, the controller asymptotically tracks the desired reference signals while maintaining the input and output constraints and minimizing steady-state control effort. The performance in Figure \ref{Fig:OSS-TwoLoop} is slightly poorer than that in Figure \ref{Fig:OSS-Inv} due to cross-coupling effects between the two dynamic loops. The robust/optimal controller of Figure \ref{Fig:OSS-TwoLoop-Optimal} does not suffer from these cross-coupling effects, and even provides an improvement over the inversion-based controller.

\subsection{Frequency Control in Power Systems}

Dynamic models of high-voltage AC power systems have the property that they are internally stable, and that the steady-state frequency deviation in the system is proportional to the net imbalance between generation and load \cite{PK:94}. This leads to a steady-state model \eqref{Eq:SteadyStateMap} of the form
\begin{equation}\label{Eq:PowerSyst}
\Delta\bar{\omega} = \tfrac{1}{\beta}\vones[m]\vones[m]^{\T}\Delta \bar{u} - \tfrac{1}{\beta}\vones[m]\bar{w}.
\end{equation}
where $\Delta \bar{\omega} \in \real^m$ denotes frequency deviations at $m$ buses of the system, $\Delta \bar{u} \in \real^m$ denotes steady-state generation change at those buses, $\bar{w}$ denotes the total loading change in the system, and $\beta > 0$ is a constant. The problem of interest is minimization of generation cost, subject to regulation of frequency, expressed as
\begin{subequations}\label{Eq:OSSFreq}
\begin{align}
\label{Eq:OSSFreq-1}
\minimize_{\Delta\bar{u},\Delta \bar{\omega}} &\quad \sum_{i=1}^{m}\nolimits J_{i}(\Delta\bar{u}_i)
\\
\label{Eq:OSSFreq-2}
\subto &\quad  0 = \beta \Delta \bar{\omega}_m
\end{align}
\end{subequations}
where $J_i$ captures the $i$th generation cost and embeds any associated unit limits. Under appropriate assumptions on $J_i$, our previous results immediately yield several provably stable controllers for the solution of this problem. The primal-dual controller \eqref{Eq:PrimalDual} reduces to
\[
\begin{aligned}
\tau_{\rm p}\Delta\dot{u}_i &= -\nabla J_i(\Delta u_i) - \mu, \qquad \tau_{\rm d}\dot{\mu} = \Delta \omega_m\\
\end{aligned}
\]
which consists of a single frequency integrator and a decentralized update for $u_i$; the inversion-based controller \eqref{Eq:Inversion} reduces to
\[
\tau\dot{\mu} = \Delta \omega_{m}, \qquad \Delta u_i = \nabla J_i^*(\mu),
\]
which precisely recovers the control scheme proposed in \cite{FD-SG:17}. For the two-loop design of Section \ref{Sec:TwoLoop}, we follow \eqref{Eq:TEquation} and must compute
\[
\mathrm{null}\begin{bmatrix}
I_r & -G_u\\
H_z & H_{u}
\end{bmatrix} = \mathrm{null}\begin{bmatrix}
I_m & -\tfrac{1}{\beta}\vones[m]\vones[m]^{\T}\\
\beta e_m^{\T} & 0
\end{bmatrix},
\]
where $e_m$ is the $m$th unit vector of $\real^m$. This nullspace is spanned by vectors of the form $\mathrm{col}(0,\xi)$ where $\xi^{\T}\vones[m] = 0$. Let $L = L^{\T}$ denote the Laplacian matrix of an undirected, connected, and weighted graph over $m$ nodes \cite[Chp. 6--8]{FB-LNS}.  In this case $\mathrm{null}(L) = \vones[m]$, and therefore $\mathrm{range}(L^{\T}) = \setdef{\xi}{\xi^{\T}\vones[m] = 0}$. If we block partition $L$ as
\[
L = \begin{bmatrix}
L_{11} & L_{12}\\
L_{21} & L_{22}
\end{bmatrix}
\]
with $L_{11} \in \real^{(m-1)\times(m-1)}$, then it follows that an eligible selection of $T$ is $T_z = 0$ and $T_u = \left[\begin{smallmatrix}
L_{11}^{\T} \\ L_{12}^{\T}
\end{smallmatrix}\right]$. Since $H_zG_u + H_u = \vones[m]^{\T}$, to satisfy Theorem \ref{Thm:TwoLoop} we may select $K_2 = e_m$. In this case, $\Pi_{\rm c} = I_m - e_m\vones[m]^{\T}$, and one may then verify that an eligible selection of $K_1$ satisfying $\Pi_{\rm c}K_1 = T_u P$ is $K_1 = \left[\begin{smallmatrix}L_{11}^{\T} \\ 0\end{smallmatrix}\right]P$. Selecting $P = L_{11}^{-1} \succ 0$, we obtain the simpler choice $K_1 = \left[\begin{smallmatrix}I \\ 0\end{smallmatrix}\right]$. With these constructions, the stabilizing controller of Theorem \ref{Thm:TwoLoop} becomes
\[
\begin{aligned}
\tau_1 \dot{\eta}_{1,i} &= -\sum_{j=1}^{m}\nolimits  a_{ij}(\nabla J_i(\Delta u_i)-\nabla J_j(\Delta u_j))\\
\tau_2 \dot{\eta}_2 &= -\Delta \omega_m\\
\Delta u &= \mathrm{col}(\eta_1,\eta_2)
\end{aligned}
\]
for $i \in \{1,\ldots,m-1\}$, where $J(\Delta u) = \sum_{i=1}^{m}J_i(\Delta u_i)$. This novel controller has the interpretation that the $m$th generation unit implements integral control on the frequency deviation, while generators $\{1,\ldots,m-1\}$ implement a distributed averaging scheme to reach agreement on their marginal costs $\nabla J_i(\Delta u_i)$ of power production. Our results herein establish that the above controllers -- and many others, which can be obtained by varying the above constructions -- lead to provably stable closed-loop systems. These stability guarantees can even be pushed to nonlinear power system models by mirroring the arguments in, e.g., \cite{FD-SG:17,JWSP:20b}.

\section{Conclusions}
\label{Section: Conclusions}

Several low-gain controller designs have been presented which can be used to drive a stable LTI system towards the solution of a linearly-constrained convex optimization problem. This renders the general OSS control framework presented in \cite{LSPL-JWSP-EM:18l} constructive for a practically-relevant class of systems, and examples have been presented to illustrate the design procedure and the flexibility inherent within it. Open directions include the extension of this framework to nonlinear and discrete-time systems, along with continued exploration of applications for optimal steady-state control.


\renewcommand{\baselinestretch}{1}
\bibliographystyle{IEEEtran}

\bibliography{brevalias,Main,JWSP,New}

\end{document}